\documentclass[12pt,reqno, oneside]{amsart}
\usepackage[text={6.5in,8.51in},centering]{geometry}
\advance\headsep by 10pt
\advance\topmargin by -10pt
\usepackage{graphicx}
\usepackage{hyperref}
\usepackage{amssymb}
\renewcommand{\S}{\mathfrak{S}}
\newcommand{\Av}{A^\bullet}
\newcommand{\Ae}{A^-}
\newcommand{\Ave}{A^{\bullet-}}
\newcommand{\Uv}{U^\bullet}
\newcommand{\Ue}{U^-}
\newcommand{\Uve}{U^{\bullet-}}

\newcommand{\C}{C}
\renewcommand{\Delta}{D}
\newcommand{\g}{G}
\newcommand{\m}{m}

\newtheorem{thm}{Theorem}[section]

\DeclareMathOperator{\fix}{fix}

\newcommand{\wF}{\widetilde{F}}
\newcommand{\wR}{\widetilde{R}}

\newcommand{\seqcite}[1]{\cite[\href{http://oeis.org/#1}{{#1}}]{oeis}}
\renewcommand{\l}{\lambda}
\renewcommand{\a}{\alpha}

\newcommand{\ipe}[2]{E_{#1}\{#2\}}

\newcommand{\bkm}[1]{Billey, Konvalinka, and Matsen \cite[#1]{tanglegrams}}

\theoremstyle{definition}

\newcommand{\qdash}[1]{\vdash_{\kern -1.6pt #1}}

\usepackage{mathtools}
\DeclarePairedDelimiter{\pbrac}{(}{)}
\DeclarePairedDelimiter{\sbrac}{[}{]}
\DeclarePairedDelimiter{\cbrac}{\{}{\}}

\newcommand*{\abs}[1]{\ensuremath{\left\vert {#1} \right\vert}}

\begin{document}
\title{\leavevmode\vspace{-10pt}Counting tanglegrams with species}
\author{Ira M. Gessel$^*$}
\address{Department of Mathematics\\
   Brandeis University\\
   Waltham, MA 02453-2700}
\email{gessel@brandeis.edu}
\date{June 22, 2021}
\thanks{$^*$Supported by a grant from the Simons Foundation (\#427060, Ira Gessel)}

\begin{abstract}
A tanglegram is a pair of binary trees with the same set of leaves. Unlabeled tanglegrams were counted recently by Billey, Konvalinka, and Matsen, who also proposed the problem of counting several variations of unlabeled tanglegrams. We use the theory of combinatorial species to solve these problems.
\end{abstract}

\maketitle

\bigskip

\section{Introduction}
A \emph{tanglegram} is a diagram, used in biology to compare phylogenetic trees, consisting of two (usually binary) trees together with a matching of their leaves.
Tanglegrams were recently counted by Billey, Konvalinka, and Matsen \cite{tanglegrams}, and we refer to this paper (and their related paper \cite{tanglegrams-tool}) for references  to biological applications. 
We answer here several questions raised by Billey, Konvalinka, and Matsen, by giving formulas for counting three variations of tanglegrams.

\begin{figure}[htbp] 
   \centering
   \includegraphics[width=1.5in]{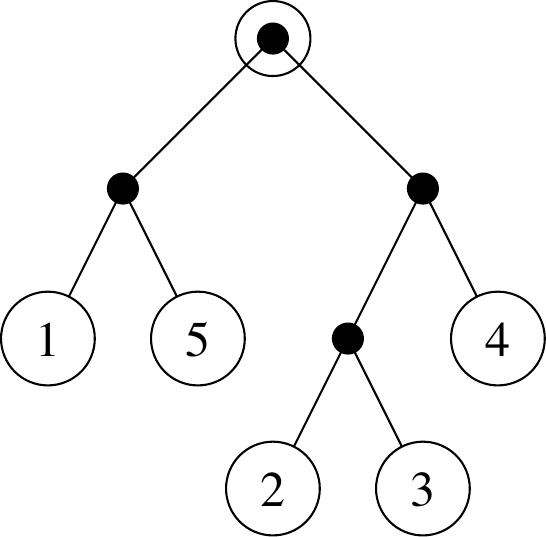} 
   \caption{A binary tree}
   \label{fig:bintree}
\end{figure}

We define a \emph{binary tree} to be a rooted tree in which every vertex has either zero or two children, and in which the leaves (vertices with no children) are labeled with distinct labels but the interior vertices are unlabeled. (See Figure \ref{fig:bintree}.)
Note that the tree with only one (labeled) vertex is a binary tree. The children of an interior vertex are not ordered, so, for example, there is one binary tree with label set $\{1,2\}$.
It is not hard to show that the number of binary trees with label set $[n]=\{1,2,\dots, n\}$ is $1\cdot3\cdots(2n-3)$ for $n>1$.  (See, e.g., \cite[Example 5.2.6]{ec2}.)

\begin{figure}[htbp] 
   \centering
   \includegraphics[width=2.5in]{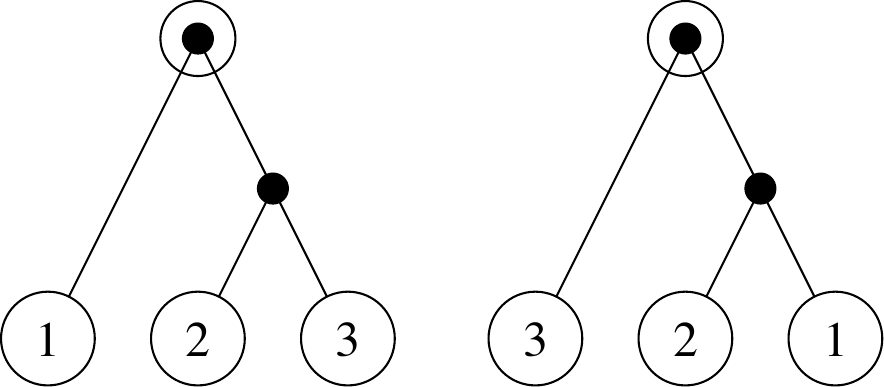} 
   \caption{A labeled tanglegram with three leaves}
   \label{fig:tanglegram1}
\end{figure}

\begin{figure}[htbp]
   \centering
   \includegraphics[width=2in]{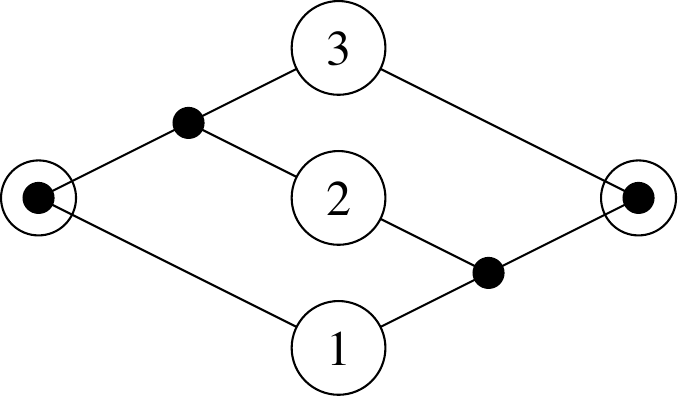} 
   \caption{Another representation of the tanglegram of Figure \ref{fig:tanglegram1}}
   \label{fig:tanglegram}
\end{figure}

We define a \emph{labeled tanglegram} to be an ordered pair of binary trees with the same set of leaf labels. Figure \ref{fig:tanglegram1} shows a labeled tanglegram with three leaves and Figure \ref{fig:tanglegram} shows another way of drawing the same tanglegram.
Labeled tanglegrams are easy to count: the number of labeled tanglegrams with $n$ leaves is $(1\cdot3\cdots(2n-3))^2$.

\begin{figure}[htbp] 
   \centering
   \includegraphics[width=2in]{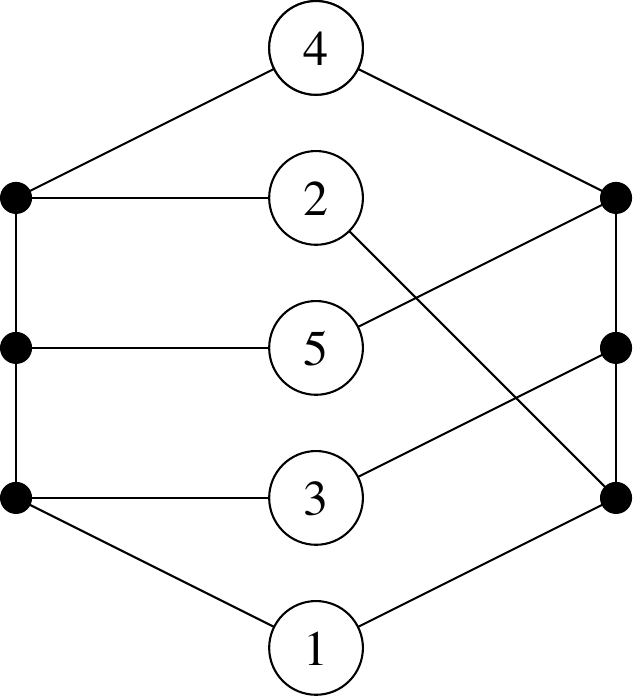} 
   \caption{An unrooted tanglegram with five leaves}
   \label{fig:unrooted-tang}
\end{figure}

An \emph{unlabeled tanglegram} is an isomorphism class of tanglegrams, where two tanglegrams are considered to be isomorphic if one can be obtained from the other by permutation of the labels.
Billey, Konvalinka, and Matsen~\cite{tanglegrams} proved a formula for the number of unlabeled tanglegrams with $n$ leaves. They left open the problem of counting several variations of unlabeled tanglegrams: unordered tanglegrams, unrooted tanglegrams, and unordered unrooted tanglegrams.  An \emph{unordered tanglegram} is an \emph{unordered} pair of binary trees (not necessarily distinct) with the same set of leaf labels and an \emph{unrooted tanglegram} is an ordered pair of unrooted trees with the same set of leaf labels in which every vertex of each tree has degree one or three, and only the leaves are labeled. (For unrooted trees, leaves are vertices of degree one.) An unrooted tanglegram is shown in  Figure \ref{fig:unrooted-tang}.

In this paper we count these three variations of tanglegrams using the theory of combinatorial species, which is an effective tool for counting unlabeled graphs of various kinds. In Section \ref{s-species} we briefly review the theory of combinatorial species; in particular we discuss the operations of Cartesian product and composition of species, and introduce two important generating functions associated to a species: the cycle index series and the unlabeled generating function. In Section \ref{s-tanglegrams}
we describe the species of rooted trees and the species of tanglegrams, which is the Cartesian product the species of rooted trees with itself. In Section \ref{s-unordered} we review the little-known operation of inner plethysm for species and symmetric functions, and apply it to the enumeration of unordered tanglegrams. Unrooted tanglegrams, both ordered and unordered, can easily be counted once we know the cycle index series for ``unrooted binary trees," which is computed in Section \ref{s-unrooted} using Leroux's dissymmetry theorem. Finally in Section \ref{s-bkm} we discuss the remarkable explicit formula found by \bkm{Proposition 4} for the cycle index series for binary trees, and we show how it can be derived by iterative application of the binomial theorem.

\section{Species}
\label{s-species}
The theory of combinatorial species, initiated by Andr\'e Joyal
\cite{joyal,joyal2},
allows us to construct combinatorial objects in ways that enable us to count them. We give here a very brief account of part of the theory; we refer the reader to Bergeron, Labelle, and Leroux \cite{species-book} for a comprehensive exposition. A concise introduction to the theory of species can be found in  \cite{Hardt}.

 A \emph{species} is a functor from the category of finite sets with bijections to itself. A species $F$ associates to  each finite set $A$ a finite set $F[A]$, called the set of \emph{$F$-structures} on $A$, and associates to  each bijection of finite sets $\sigma\colon A \to B$ a bijection $F[\sigma]\colon F[A]\to
 F[B]$. 
In particular, any bijection $\pi\colon[n]\to[n]$ yields a bijection $F[\pi]\colon F[n]\to F[n]$, where $F[n] =F[\{ 1,2,\dots,n\}]$, so  the symmetric group $\mathfrak{S}_n$ acts on
the set $F[n]$.
The $\mathfrak{S}_n$-orbits under
this action are called \emph{unlabeled}
$F$-structures of order $n$.
To any species $F$ we may associate its 
\emph{cycle index series} $Z_F$, a symmetric function
defined by
\begin{equation}
\label{eq:cycinddef}
  Z_F = Z_{F}(p_{1}, p_{2}, \dots) = \sum_{n=0}^\infty\frac{1}{n!} \biggl(\, \sum_{\sigma \in \S_{n}} 
     \fix  {F \sbrac{\sigma}}\, p_{\sigma} \biggr).
 \end{equation}
 Here $p_k$ is the \emph{power sum symmetric function%
\footnote{In many accounts of cycle indices our power sums $p_i$ are replaced with  
 indeterminates. In particular, our $p_i$ is written as $x_i$ in  \cite{species-book}.
\endgraf
The sum of the terms of degree $n$ in $Z_F$ form the ``characteristic"   
 or ``Frobenius image" of the representation of $\S_n$ associated with the action of 
 $\S_n$ on $F[n]$ (see, e.g., \cite[pp.~351--352 and 395--396]{ec2}) and all of the 
 operations on species that we discuss have analogues for representations of symmetric groups.  This is one of the reasons why we consider cycle indices to be symmetric functions.
}}
$\sum_{j=1}^\infty x_j^k$,
$\fix {F \sbrac{\sigma}} = \abs{\cbrac{s \in F \sbrac{n} : F \sbrac{\sigma} \pbrac{s} = s}}$,
and $p_{\sigma} = p_{1}^{\sigma_{1}} p_{2}^{\sigma_{2}} \dots$, where $\sigma_i$
is the number of $i$-cycles of $\sigma$.
 
The inner sum in \eqref{eq:cycinddef} contains $n!$ terms, many of which are equal. It is much more efficient to replace this sum with a sum over partitions of $n$. To do this, we associate to
each $\sigma\in\mathfrak{S}_n$ a partition of $n$, called its \emph{cycle type}, by listing the lengths of its cycles  in weakly decreasing order. Thus the cycle type of the permutation $(142)(5)(3)$ is $(3,1,1)$. Two permutations in $\mathfrak{S}_n$ are conjugate if and only if they have the same cycle type. It is not hard to show that $\fix F[\sigma]$ depends only on the conjugacy class of $\sigma$, so for each partition $\l$ of $n$ we may define $\fix F[\l]$ to be $\fix F[\sigma]$ where $\sigma$ is a permutation of $[n]$ with cycle type $\l$. We define the power sum symmetric functions $p_\l$ similarly, so if $\lambda=(\lambda_1, \lambda_2,\dots, \lambda_k)$ then  $p_\l=p_{\l_1}p_{\l_2}\cdots p_{\l_k}$.
If the partition $\lambda$ has $m_i$ parts equal to $i$, for each $i$, then we define $z_\l$ to be $z_\l=1^{m_1} m_1!\, 2^{m_2} m_2!\cdots$, so the number of partitions in $\mathfrak{S}_n$ of cycle type $\l$ is $n!/z_\l$.
Thus a more compact formula for the cycle index $Z_\l$ is 
\[Z_F=\sum_{n=0}^\infty
 \biggl(\sum_{\l \vdash n} \fix F[\l]\, \frac{p_{\l}}
 {z_{\l}}
 \biggr),\]
where  $\l\vdash n$ means that the sum is over all partitions $\l$ of $n$. 

For any species $F$, we denote by $\wF(x)$ the ordinary generating function for unlabeled $F$-structures; that is,
\begin{equation*}
\wF(x) = \sum_{n=0}^\infty f_n x^n,
\end{equation*}
where $f_n$ is the number of unlabeled $F$-structures of order $n$. Then
by \cite[p.~18, Theorem 8b]{species-book} we have
\begin{align}
\wF(x)
&=Z_F(x, x^2, x^3, \dots)\label{e-unl1}\\
&=\sum_{n=0}^\infty x^n
 \biggl(\sum_{\l \vdash n} \frac{\fix F[\l]}
 {z_{\l}}
 \biggr).\notag
\end{align}
That the number of unlabeled $F$-structures of order $n$ is 
$
\sum_{\l \vdash n} {\fix F[\l]}/
 {z_{\l}}
$
can also be seen directly by Burnside's lemma (see Section \ref{s-tanglegrams}).

Among the most important species are the species $E_n$ of $n$-element sets, defined 
for $n\in \mathbb{N}$ by
\begin{equation*}
E_n[A]=\begin{cases} 
\{A\},&\text{if $|A|=n$}\\
\varnothing,&\text{otherwise}.
\end{cases}
\end{equation*}
The cycle index of $E_n$ is the \emph{complete symmetric function} $h_n$, defined by 
\begin{equation*}
h_n = \sum_{i_1\le i_2\le \dots\le i_n}x_{i_1}x_{i_2}\cdots x_{i_n},
\end{equation*}
and also given by the formula
\begin{equation}
\label{e-ph}
h_n = \sum_{\l\vdash n} \frac{p_\l}{z_\l},
\end{equation}
which follows from 
\begin{equation*}
\sum_{n=0}^\infty h_n = \prod_{j=1}^\infty \frac{1}{1-x_j} =\exp\biggl(\sum_{k=1}^\infty \frac{p_k}{k}\biggr).
\end{equation*}

The special case $E_1$, the species of \emph{singleton sets}, is denoted by $X$. It has cycle index $Z_X = p_1$.

Given species $F$ and $G$, we can combine them to get the \emph{sum} $F+G$, the \emph{product} $F  G$,  the \emph{composition} $F(G)$ (also denoted $F\circ G$), and the 
\emph{Cartesian product} $F\times G$, and these operations on species translate into operations on cycle indices. We refer the reader to \cite[pp.~1--58]{species-book} for details about these operations.

The sum $(F+G)[A]$ is the disjoint union of $F[A]$ and $G[A]$. 
An $FG$-structure on the set $A$ is obtained by partitioning $A$ into disjoint subsets $B$ and $C$ (possibly empty) and taking an $F$-structure on $B$ and a $G$-structure on $C$. 

Suppose that%
\footnote{The composition $F(G)$ can be defined without this condition, but we will not need the more general definition.}
$G[\varnothing]=\varnothing$.
Then an  $F(G)$-structure on the set $A$ is an $F$-structure of $G$-structures; more precisely, $F(G)[A]$ is the set of triples  $(\pi, \alpha, \beta)$, where
$\pi$ is a partition of the set $A$, $\alpha$ is an $F$-structure on (the set of blocks of) $\pi$,
and $\beta$ is a set of $G$-structures on the blocks of $\pi$. In particular, an $E_n(G)$-structure on $A$ is a partition of $A$ into $n$ blocks, together with a $G$-structure on each block.

The Cartesian product is defined by $(F\times G)[A]=F[A]\times G[A]$; thus an $F\times G$-structure is a pair of structures on the same set.

The corresponding operations for cycle indices are simple for the sum and product: 
$Z_{F+G}=Z_{F}+Z_{G}$ and $Z_{FG}=Z_{F}Z_{G}$. The cycle index operation for composition of species is an operation on symmetric functions called \emph{composition} or \emph{plethysm} (see, e.g., \cite[p.~447]{ec2}). The composition of symmetric functions $f$ and $g$,  denoted by $f[g]$ or  $f\circ g$, may be defined by 
\begin{equation*}
f[g] = f(g(p_1, p_2, p_3, \dots, ), g(p_2, p_4, p_6,\dots), \dots);
\end{equation*}
i.e., $f[g]$ is obtained from $f$ by  replacing each $p_i$ with $p_i[g]=g(p_i, p_{2i}, p_{3i},\dots)$.
Then $Z_{F(G)}=Z_F[Z_G]$.
In particular, since $Z_{E_2}=h_2 = \tfrac12(p_1^2+p_2)$, we have 
$Z_{E_2}[g] = \tfrac12(g^2 + p_2[g])$. The cycle index operation corresponding to the Cartesian product on species is an operation on symmetric functions called the \emph{Kronecker product}, \emph{internal product}, or \emph{inner product}.
The Kronecker product, denoted by $*$, is defined by 
$p_\l*p_{\mu}=\delta_{\mu\l}z_{\l}p_\l$ and linearity, or equivalently,
\begin{equation*}
\sum_{\l} a_{\l} \frac{p_\l}{z_\l} *
\sum_{\l} b_{\l} \frac{p_\l}{z_\l}
=\sum_{\l} a_{\l}b_{\l}\frac{p_\l}{z_\l}.
\end{equation*}
Then $Z_{F\times G}=Z_F * Z_G$.

As is customary in discussing species we will consider isomorphic species to be equal; for example, in equation \eqref{e-R1} below the two sides are really isomorphic rather than equal.

\section{Tanglegrams}
\label{s-tanglegrams}
Let $R$ be the species of (rooted) binary trees with labeled leaves and unlabeled internal vertices. A binary tree is either a single labeled vertex or an unlabeled root together with an unordered pair of binary trees. Thus $R$ satisfies the equation 
\begin{equation}
\label{e-R1}
R = X + E_2(R),
\end{equation}
so the cycle index $Z_R$ satisfies
\begin{equation}
\label{e-ZR}
Z_R= p_1 + h_2[Z_R].
\end{equation}
Terms of $Z_R$ can be computed fairly easily by successive substitution in \eqref{e-ZR}, though there are other ways to compute them that are  more efficient. (See, e.g., \cite[Corollary D1]{MR927766} and Section \ref{s-bkm}.) The first few terms of $Z_R$ are 
\begin{equation*}
p_{{1}}+ \left( \frac12 p_{{1}} ^{2}+\frac12 p_{{2}} \right) +
 \left( \frac12 p_1 p_{{2}}+\frac12 p_{1}^{3} \right) +
 \left( \frac58 p_{{1}}^{4}+\frac38 p_{{2}}^{2}
 +\frac34 p_{{1}}^{2}p_{{2}}+\frac14 p_{{4}} \right) +\cdots
\end{equation*}
It is easy to see from \eqref{e-ZR} that if the power sum $p_i$  occurs in $Z_R$,  then $i$ is a power of~2. From \eqref{e-ZR} we can also easily derive the well-known functional equation for the ordinary generating function $\wR(x)$ \begin{equation*}
\wR(x) = x +\frac12\bigl(\wR(x)^2 + \wR(x^2)\bigr),
\end{equation*}
(see \seqcite{A001190})
but to count unlabeled tanglegrams we need the full cycle index.

Now let $T$ be the species of (labeled) tanglegrams. Since a  tanglegram   is a pair of binary trees with the same set of labels, $T$ is the Cartesian product $R\times R$, so $Z_T = Z_{R\times R}=Z_R*Z_R$. Thus if 
\begin{equation}
\label{e-r1}
Z_R = \sum_{\l}r_\l \frac{p_\l}{z_\l},
\end{equation}
where the sum is over all partitions $\l$,
then $Z_T = \sum_\l r_\l^2 p_\l/z_\l$, and the number of unlabeled tanglegrams with $n$ leaves is
\begin{equation}
\label{e-ut1}
\sum_{\l\vdash n} \frac{r_\l^2}{z_\l}.
\end{equation}
A formula equivalent to \eqref{e-ut1} was given by \bkm{Theorem 1}; we will discuss their result further in Section \ref{s-bkm}.

A \emph{tangled chain of length $k$}   is a $k$-tuple of binary trees sharing the same set of leaves. It is clear that the species of tangled chains is the $k$th Cartesian power of $R$, 
so the number of unlabeled tangled chains of length $k$ is
\begin{equation}
\label{e-tc}
\sum_{\l\vdash n} \frac{r_\l^k}{z_\l},
\end{equation}
as also shown by \bkm{Theorem 3}.

It may be noted that  \eqref{e-tc} is an easy consequences of Burnside's lemma: if a group $G$ acts on a finite set $S$ then the number of orbits is
\begin{equation}
\label{e-burnside}
\frac{1}{|G|}\sum_{g\in G}\fix g,
\end{equation}
where $\fix g$ is the number of elements of $S$ fixed by $g$.
To derive \eqref{e-tc} from Burnside's lemma, we consider the action of $\S_n$ on $k$-tuples of labeled binary trees with leaf set $[n]$. A $k$-tuple is fixed by a permutation if and only if all its entries are fixed, so the $n!/z_\lambda$ permutations of cycle type $\lambda$ contribute $(n!/z_\lambda)r_\lambda^k$ to the sum \eqref{e-burnside}. Thus the number of orbits is 
\begin{equation*}
\frac{1}{n!}\sum_{\lambda\vdash n}\frac{n!}{z_\lambda} r_\lambda^k
=\sum_{\lambda\vdash n}\frac{r_\lambda^k}{z_\lambda}.
\end{equation*} 

\section{Unordered tanglegrams}
\label{s-unordered}

To count unordered tanglegrams, we use another operation on species, \emph{inner plethysm}, that is not as well known as the other operations. Inner plethysm is a kind of composition of species that bears the same relation to the Cartesian product that ordinary composition bears to the ordinary product. It is closely related to the operation of functorial composition of species introduced in \cite{functorial} and discussed further in \cite[Section 2.2]{species-book}. The term ``inner plethysm" was introduced by D.~E.~Littlewood \cite{MR0095208} for the corresponding operation on symmetric functions, and the species operation was introduced by L.~Travis in his Ph.$\,$D. thesis \cite{MR2698697}.  We refer to this thesis for results about inner plethysm not proved here. 

There is no  standard notation for inner plethysm, so we will introduce the notation $F\{G\}$ for the inner plethysm of species, with the same notation for inner plethysm of symmetric functions. We will define here only the  inner plethysm $\ipe nG$, which is all that we  need: for any finite set $A$, $\ipe nG[A]$ is the set of multisets of size $n$ of elements of $G[A]$. We can define $\ipe nG$ in another way: The symmetric group $\S_n$ acts on the elements of the $n$th Cartesian power $G^{\times n}[A]$ by permuting the $n$ entries, and the elements of $\ipe nG[A]$ are the orbits under this action. (The functorial  composition of species is defined similarly, but with a \emph{set}, rather than a multiset, of elements of $G[A]$.)

Inner plethysm of symmetric functions is determined by the following:
\begin{enumerate}
\item for fixed $g$, the map $f\mapsto f\{g\}$ is a homomorphism from the ring of symmetric functions with the usual product to the ring of symmetric functions with the Kronecker product
\item For a partition $\l$ and an integer $k$, let $\l^k$ denote the cycle type of the $k$th power of a permutation with cycle type $\l$. Then
\[
p_k\biggl\{\sum_{\l}a_{\l}\frac{p_{\l}}{z_\l}\biggr\}
  =\sum_{\l}a_{\l^k}\frac{p_{\l}}{z_\l}.
\]
\end{enumerate}
Travis \cite[Theorem 2.12]{MR2698697} showed that for any species $F$ and $G$, we have $Z_{F\{G\}}=Z_F\{Z_G\}$. 

It is clear that the species  of unordered tanglegrams is $\ipe 2R$, where $R$ is the species of binary trees. So the cycle index for unordered tanglegrams is 
$h_2\{Z_R\}=\tfrac12(p_1^2+p_2)\{Z_R\}$.
Thus if $Z_R=\sum_{\l}r_\l p_\l/z_\l$ then the cycle index for unordered tanglegrams is 
\begin{equation*}
\frac12\biggl(\sum_\l r_\l^2\frac{p_\l}{z_\l}\biggr)
+\frac12\biggl(\sum_\l r_{\l^2}\frac{p_\l}{z_\l}\biggr),
\end{equation*}
and we obtain the number of unordered tanglegrams with $n$ leaves by setting each $p_\l$ to~1 in the sum of the terms of degree $n$. 
So if $a_n$ is the number of unordered tanglegrams with $n$ leaves then
\begin{equation*}
a_n = \frac12 \sum_{\l \vdash n}\frac{r_\l^2 + r_{\l^2}}{z_\l}.
\end{equation*}

This formula for unordered tanglegrams can also be derived directly from Burnside's lemma, using the action of $\S_n\times \S_2$ on labeled tanglegrams, where $\S_2$ acts by permuting the two trees.

Here are the first few values of the number $a_n$ of unordered tanglegrams with $n$ leaves (see \seqcite{A259114}):
\[\vbox{\halign{\ \hfil\strut$#$\hfil\ \vrule&&\hfil\ $#$\ \hfil\cr
n&1&2&3&4&5&6&7&8&9&10&11\cr
\noalign{\hrule}
a_n&1&1&2&10&69&807&13048&269221&6660455&191411477&6257905519\cr}}\]

Similarly, $\ipe kR$ is the species of unordered tangled chains of length $k$.

\section{Unrooted tanglegrams}
\label{s-unrooted}
To count unrooted tanglegrams, we need to find the cycle index for unrooted trees, and to do this we  use a \emph{dissymmetry theorem}. Dissymetry theorems, introduced by Leroux \cite{leroux} reduce the enumeration of unrooted trees to the enumeration of several types of rooted trees, which can usually be counted through decompositions. 
The basic dissymmetry theorem says that
if $A$ is a species of unrooted trees of some type, $\Av$ is the species of  $A$-trees rooted at a vertex%
\footnote{A tree rooted at a vertex is formally an ordered pair $(T,v)$, where $T$ is a tree and $v$ is a vertex of $T$. Trees rooted at edges are defined similarly.}, $\Ae$ is the species of $A$-trees rooted at an edge, and $\Ave$ is the species of $A$-trees rooted at a vertex and incident edge (or equivalently, at a directed edge), then 
\begin{equation}
\label{e-diss}
A+\Ave=\Av+\Ae.
\end{equation}
We give here a brief sketch of the proof of \eqref{e-diss}, referring the reader to \cite[Section 4.1]{species-book} for a more detailed discussion. Every tree has a unique ``center", which is  a vertex or edge that is fixed by every automorphism of the tree. An unrooted tree may be identified with a tree rooted at its center. To prove \eqref{e-diss}, we describe a bijection, equivariant with respect to the automorphism group of the tree, from the non-center vertices and edges of a tree to pairs consisting of a vertex and an incident edge. If $v$ is a non-center vertex, we pair it with the first edge on the unique path from $v$ to the center (this edge may be the center), and if $e$ is a non-center edge, we pair it with the first vertex on the unique path from $e$ to the center (this vertex may be the center). From the bijection just described, we get a bijection from  $A$-trees rooted at a vertex or edge to $A$-trees rooted at a center vertex or edge (equivalent to unrooted $A$-trees) or at a vertex and incident edge.

Now let $U$ be the species of ``unrooted binary trees"; that is, unrooted trees in which every vertex has degree one or three, the leaves (vertices of degree one) are labeled, and the internal vertices (of degree three) are unlabeled. (We are not including the tree with one vertex.)  
First we consider $U$-trees rooted at an edge $e$. Removing the  edge $e$ and rooting 
the remaining two trees at the vertices incident with $e$ gives two rooted binary trees. Thus $\Ue=E_2(R) = R-X$, where $R$ is the species of rooted binary trees discussed in Section \ref{s-tanglegrams}.
 Similarly, we can remove the root edge from a $U$-tree rooted at a vertex and incident edge to obtain a pair of rooted trees, but in this case the pair is ordered, so $\Uve=R^2$. Finally, the $U$-trees rooted at a vertex may be rooted at either an internal vertex or a leaf. The species of $U$-trees rooted at an internal vertex is $E_3(R)$ and the species of $U$-trees rooted at a leaf is $XR $, so $\Uv = E_3(R)+XR$.
Thus \eqref{e-diss} gives 
\begin{equation*}
U+R^2 = E_3(R) +XR +R-X,
\end{equation*}
and we obtain a formula for the cycle index of $U$,
\begin{equation*}
Z_U =h_3[Z_R]+p_1 Z_R+Z_R -Z_R^2- p_1.
\end{equation*}
The first few terms of $Z_U$ are
\begin{equation*}
\left( \frac12 p_{1}^{2}+\frac12p_{{2}} \right)+
 \left(\frac16p_{{1}}^{3}+ \frac12 p_1 p_{{2}}+\frac13p_{{3}} \right) 
 + \left(\frac18p_{{1}}^{4}+ \frac14 p_{{1}}^{2} p_{{2}}+\frac38p_{{2
}}^{2}+\frac14p_{{4}} \right)+\cdots
\end{equation*}

Then the species of unrooted tanglegrams is the Cartesian product $U\times U$, with cycle index $Z_U* Z_U$ and the species of unrooted unordered tanglegrams is $E_2\{U\}$, with cycle index $h_2\{Z_U\}$.
These cycle indices are easily computed, and the  numbers of unlabeled tanglegrams of these types may be obtained from the cycle indices by \eqref{e-unl1}.
The numbers $b_n$ of unrooted tanglegrams and $c_n$ of unrooted unordered tanglegrams, with $n$ leaves, for small values of $n$ are as follows:

\[\vbox{\halign{\ \hfil\strut$#$\hfil\ \vrule&&\hfil\ $#$\ \hfil\cr
n&2&3&4&5&6&7&8&9&10&11&12\cr
\noalign{\hrule}
b_n&1&1&2&4&31&243&3532&62810&1390718&36080361&1076477512\cr}}\]

\vspace{-8pt}

\[\vbox{\halign{\ \hfil\strut$#$\hfil\ \vrule&&\hfil\ $#$\ \hfil\cr
n&2&3&4&5&6&7&8&9&10&11&12\cr
\noalign{\hrule}
c_n&1&1&2&4&22&145&1875&31929&698183&18056523&538340256\cr}}\]
The sequence $b_n$ is \seqcite{A259115} and the sequence $c_n$ is \seqcite{A259116}.

Similarly, the $k$th Cartesian power $U^{\times k}$ is the species of unrooted tangled chains of length $k$ and the inner plethysm $\ipe kR$ is the species of unordered unrooted tangled chains of length $k$.

Symmetric function computations were done with the help of John Stembridge's Maple package for symmetric functions \cite{SF,SF2.4}.

\section{The formula of Billey, Konvalinka, and Matsen}
\label{s-bkm}
\bkm{Proposition 4} proved (though they did not state it this way) that if $\l$ is a binary partition of $n$, that is, a partition in which every part is a power of 2, then
the coefficient $r_\l$ of $p_\l/z_\l$ in $Z_R$, as defined in \eqref{e-ZR}, is given by the simple explicit formula
\begin{equation}
\label{e-bkmr}
r_\l=\prod_{i=2}^{l(\l)}\bigl(2(\l_i+\cdots + \l_{l(\l)}) -1\bigr),
\end{equation}
where $l(\l)$ is the number of parts of $\l$. (If $\l$ is not a binary partition then $r_\l=0$.) They proved \eqref{e-bkmr} by showing that the right side of \eqref{e-bkmr} satisfies the same recurrence that \eqref{e-ZR} implies for  the coefficients of $Z_R$. A bijective proof of \eqref{e-bkmr} was given by Fusy \cite{fusy}.

Although  Billey, Konvalinka, and Matsen's proof of \eqref{e-bkmr} is short and direct, it is of interest to place \eqref{e-bkmr} into a broader context.

If we let $r(z)$ be the result of setting $p_1=z$ and $p_i=0$ for $i>1$ in $Z_R$ then, as noted in~\cite {tanglegrams},  $r(z)$ is the exponential generating function for labeled binary trees, and \eqref{e-ZR} yields
\begin{equation}
\label{e-rquad}
r(z) =z +r(z)^2/2.
\end{equation}
There are two ways of deriving from \eqref{e-rquad} the explicit formula 
\begin{equation}
\label{e-rzf}
r(z) 
  = \sum_{n=1}^\infty 1\cdot 3\cdot 5\cdots (2n-3) \frac{z^n}{n!},
\end{equation}
in which the coefficient of $z^n\!/n!$ is given by \eqref{e-bkmr} for $\lambda = (1^n)$; we might hope that at least one of them will generalize to give \eqref{e-bkmr}.

First, we can apply Lagrange inversion \cite{lagrange}, which can solve functional equations of the form $f(z) = z + G\bigl(f(z)\bigr)$, and we obtain
\[r(z) = \sum_{n=1}^\infty \frac{1}{2^{n-1}n}\binom{2n-2}{n-1} z^n
  = \sum_{n=1}^\infty 1\cdot 3\cdot 5\cdots (2n-3) \frac{z^n}{n!},
\]
and more generally, for any positive integer $k$ we have
\begin{align*}
r(z)^k &= \sum_{n=k}^\infty \frac{k}{2^{n-k}(2n-k)}\binom{2n-k}{n} z^n\\
  &=k!\sum_{n=k}^\infty 1\cdot3\cdot 5\cdots (2n-2k-1)\binom{2n-k-1}{k-1} \frac{z^n}{n!}.
\end{align*}
Here and elsewhere an empty product is taken to be 1.
  There is a generalization of Lagrange inversion  for plethystic equations, due to Labelle \cite{MR927766} (see also \cite{MR1354969}), which can be applied to \eqref{e-ZR}. However,  it is not clear how \eqref{e-bkmr} can be obtained from Labelle's formula.

The second approach to 
equation \eqref{e-rzf} is through the binomial theorem. We can solve \eqref{e-rquad} to get 
$r(z) =1-\sqrt{1-2z}$, and then apply the binomial theorem to get
\begin{equation*}
r(z) = \sum_{n=1}^\infty (-1)^{n-1}\binom{\tfrac12}n 2^n z^n 
  = \sum_{n=1}^\infty 1\cdot 3\cdot 5\cdots (2n-3) \frac{z^n}{n!}.
\end{equation*}

More generally, the binomial theorem gives for any $k$
\begin{equation*}
\bigl(1-r(z)\bigr)^k = (1-2z)^{k/2}=1-\sum_{n=1}^\infty k(2-k)(4-k)\cdots(2n-2-k)\frac{z^n}{n!}.
\end{equation*}

If generalized Lagrange inversion is the right way to derive \eqref{e-bkmr}, then one would expect to find similar nice formulas for $Z_R^k$, for  integers $k>1$.
On the other hand, if \eqref{e-bkmr} is a consequence of a binomial theorem for symmetric functions then we might expect to find nice formulas for coefficients of $(1-Z_R)^k$ rather than for $Z_R^k$.

It turns out that we find nice formulas for $(1-Z_R)^k$ but not for $Z_R^k$, suggesting that \eqref{e-bkmr} may  be related to the binomial theorem, and this is in fact the case.
. 

We now show how  formula \eqref{e-bkmr}, and more generally a formula for the coefficients of arbitrary powers of $1-Z_R$, can be derived from the ordinary binomial theorem, using an iterative approach similar to that used by Wagner \cite{MR2255414} in deriving a formula for the cycle index series for rooted trees.

In the following proof, any product of the form $\prod_{i=n+1}^{n}u_i$ is interpreted as 1 and any sum of the form $v_j+v_{j+1}+\cdots + v_n$ with $j=n+1$ is interpreted as 0.

\begin{thm}
\label{t-bin}
Let $\g= 1-Z_R$, where $Z_R$ satisfies \eqref{e-ZR}. Then $\g$ satisfies
\begin{equation}
\label{e-g1}
\g = (p_2[\g] - 2p_1)^{1/2},
\end{equation}
and for any $\alpha$, 
\begin{equation*}
\g^{-\a}=\sum_\l \frac{p_\l}{z_\l} \prod_{j=1}^{l(\l)} (\a+2\l_{j+1}+2\l_{j+2}+\cdots+2\l_{l(\l)}),
\end{equation*}
where the sum is over all binary partitions $\l$. 
\end{thm} 

\begin{proof}
We first note that Billey, Konvalinka, and Matsen's formula \eqref{e-bkmr} may be obtained from the theorem by rewriting the product as 
\begin{equation*}
\a\prod_{j=2}^{l(\lambda)} (\a+2\lambda_j+2\lambda_{j+1}+\cdots +2\l_{l(\lambda)}),
\end{equation*}
and setting $\a=-1$.

Since $h_2=\tfrac12(p_1^2 +p_2)$, equation \eqref{e-ZR} may be written
\begin{equation*}
1-\g = p_1 + \tfrac12(p_1^2 +p_2)[1-\g] = p_1+ \tfrac12(1-\g)^2 +\tfrac12(1-p_2[\g]).
\end{equation*}
This may be rearranged to
\begin{math}
\g^2 = p_2[\g]-2p_1,
\end{math}
which yields \eqref{e-g1}.

Keeping in mind that $\g$, and thus $p_2[\g]$, has constant term 1, we may apply the binomial theorem to \eqref{e-g1} to get 
\begin{align}
\g^{-\a}&=(p_2[\g]-2p_1)^{-\a/2}=\sum_{\m=0}^\infty (-2)^{\m}\binom{-\a/2}{\m}p_1^{\m}p_2[\g]^{-\a/2-\m}\notag\\
  &=\sum_{\m=0}^\infty \frac{\a(\a+2)(\a+4)\cdots (\a+2\m-2)}{\m!}p_1^{\m}p_2[\g]^{-\a/2-\m},
\label{e-gbinom}
\end{align}
where for $\m=0$, the coefficient in \eqref{e-gbinom} is 1. 

Now let
$\C(\l,\a)$ be the coefficient of $p_\l$ in $\g^{-\a}$. Then $\C(\l,\a)$ is a polynomial in~$\a$, and is 0 unless $\l$ is a binary partition. If $\l$ is a binary partition in which 1 occurs  $\m$ times as a part, then there is a  (possibly empty) binary partition $\mu=(\mu_1,\mu_2,\dots, \mu_{l(\l) -\m})$  such that
\begin{equation*}
\l=(2\mu_1, \cdots, 2\mu_{l(\l)-\m}, \underbrace{1,\dots, 1}_{\m}).
\end{equation*}

From \eqref{e-gbinom} we see that $\C(\l,\a)$ satisfies, and is determined by, the recurrence
\begin{equation}
\label{e-rec1}
\C(\l,\a) = \frac{\a(\a+2)(\a+4)\cdots (\a+2\m-2)}{\m!}\C(\mu,\a/2+\m),
\end{equation}
together with the initial condition $\C(\varnothing, \a) = 1$. Here the factor in front 
of $\C(\mu,\a/2+\m)$ is taken to be 1 for $\m=0$.

Define $\Delta(\l,\a)$ for a binary partition $\l$ by
\begin{equation}
\label{e-Pi}
\Delta(\l,\a)=\frac{1}{z_\l}\prod_{j=1}^{k} (\a+2\lambda_{j+1}+2\lambda_{j+2}+\cdots +2\l_{k}),
\end{equation}
where we have written $k$ for $l(\l)$. We will prove that $\C(\l,\a)=\Delta(\l,\a)$ by showing that $\Delta(\l,\a)$ satisfies the same recurrence as $\C(\l,\a)$.

Suppose that $\l$ is a nonempty binary partition with $\m$ parts equal to 1. Removing from $\l$ the parts equal to 1 and then dividing each remaining part by 2 leaves another (possibly empty) binary partition $\mu$. 
A straightforward computation shows that $z_\l = \m!\, 2^{k-\m}z_\mu$.

The product of the last $\m$ factors in \eqref{e-Pi} is 
\begin{equation*}
\prod_{j=k-\m+1}^k (\a+2\lambda_{j+1}+2\lambda_{j+2}+\cdots +2\l_{k}).
\end{equation*}
If any parts of $\lambda$  occur in this product (i.e., if $m\ge2$) then the lowest indexed part that occurs is $\lambda_{k-m+2}=1$.
Thus all of the parts of $\l$ that occur in this product are equal to 1, so the product is equal to 
\begin{equation}
\label{e-prod1}
(\a+2\m-2)(\a+2\m-4)\cdots(\a+2)\a,
\end{equation}
and this is also true if $m=0$ or $1$, where the product is 1 if $m=0$ and is $\alpha$ if $m=1$.
Thus \eqref{e-Pi} is equal to \eqref{e-prod1} times
\allowdisplaybreaks
\begin{align*}
\frac{1}{z_\l}\prod_{j=1}^{k-\m}\bigl(\a+(2\l_{j+1}+&\cdots +2\l_{k-\m})+(2\l_{k-\m+1}+\cdots+2\l_{k})\bigr)\\
 &=\frac{1}{z_\l}\prod_{j=1}^{k-\m}(\a+2\m+2\l_{j+1}+\cdots +2\l_{k-\m})\\
 &=\frac{2^{k-\m}}{z_\l}\prod_{j=1}^{k-\m}\left(\a/2+\m+2\mu_{j+1}+\cdots+2\mu_{k-\m}\right)\\
 &=\frac{1}{\m!\,z_\mu}\prod_{j=1}^{k-\m}\left(\a/2+\m+2\mu_{j+1}+\cdots+2\mu_{k-\m}\right)\\
 &= \frac{1}{\m!}\Delta(\mu, \a/2+\m).
\end{align*}
Thus $\Delta(\l,\a)$ satisfies the same recurrence and initial conditions as $\C(\l,\a)$, so 
$\C(\l,\a) = \Delta(\l,\a)$. 
\end{proof}

\bigskip
\noindent\textbf{Acknowledgments.} I would like to thank Sara Billey, Gilbert Labelle,  Erick Matsen, and two anonymous referees for helpful comments on  earlier versions of this paper.


\end{document}